\documentclass{amsart}
\usepackage{amssymb}
\usepackage{amscd}

\newtheorem{thm}{Theorem}[subsection]
\newtheorem{lem}[thm]{Lemma}
\newtheorem{cor}[thm]{Corollary}
\newtheorem{prop}[thm]{Proposition}  
\newtheorem{mthm}{Theorem 4.2.2}  
  
\newtheorem{mprop}{Proposition 4.3.7}  
\newtheorem{propy}[thm]{Property}   
\newtheorem{propyI}{Property I}  
\newtheorem{propyII}{Property II}

\theoremstyle{remark}
\newtheorem{ack}{Acknowledgments}

\theoremstyle{definition}
\newtheorem{defn}[thm]{Definition}
\newtheorem{rem}[thm]{Remark}

\newtheorem{eg}[thm]{Examples}       
\newtheorem{pf}{Proof}                 
\newtheorem{conj}[thm]{Conjecture}    
   
\newtheorem{prob}[thm]{Problem}    
\newtheorem{dsconj}{Direct Sum Conjecture}  
\newtheorem{dconj}{Dimension Conjecture}  
\newtheorem{mconj}{Conjecture 2.4.1}  
  
\newtheorem{mmmconj}{Conjecture 4.3.8}

\title[MZV algebra]{The multiple zeta value algebra \\ and  
  the stable derivation algebra}
\author{Hidekazu Furusho}

\address{Research Institute for Mathematical Sciences\\ 
Kyoto University, Kyoto 606-8502 Japan}
\email{furusho@kurims.kyoto-u.ac.jp}

\begin{document}
\bibliographystyle{amsalpha+}
\maketitle

\begin{abstract}
The MZV algebra is 
the graded algebra over ${\bold Q}$ generated by all multiple zeta values.
The stable derivation algebra is a graded Lie algebra version 
of the Grothendieck-Teichm\"{u}ller group.
We shall show that there is  a canonical surjective $\bold Q $-linear map
from the graded dual vector space of the stable derivation algebra
over $\bold Q$
to the new-zeta space, 
the quotient space of the sub-vector space of the MZV algebra
whose grade is greater than $2$
by the square of the maximal ideal.
As a corollary, we get an upper-bound for the dimension of the graded
piece of the MZV algebra at each weight
in terms of the corresponding dimension of the graded piece of the 
stable derivation algebra.
If some standard conjectures by Y. Ihara and P. Deligne 
concerning the structure of the
stable derivation algebra hold, 
this will become a bound conjectured in Zagier's talk at 1st 
European Congress of Mathematics.
Via the stable derivation algebra,
we can compare the new-zeta space with the
$l$-adic Galois image Lie algebra
which is associated with
so the Galois representation on the pro-$l$ fundamental group of
$\bold P ^1_ {\overline{\bold Q}}-\left\{0,1,\infty\right\}$ .
\end{abstract}

\tableofcontents

\setcounter{section}{-1}
\section{Introduction}

In connection with the Galois representation on the pro-$l$ fundamental group
$\pi_1^l(\bold P ^1_ {\overline{\bold Q}}- \left\{0,1,\infty\right\})$
for each prime $l$, 
the {\it stable derivation algebra}
$\Bigl
(\frak D\centerdot=\underset{w> 2}{\oplus}\frak D_w ,\left\{\frak F^m \frak D\centerdot\right\}_{m\in \bold N}
\Bigr)$ 
which is a certain filtered graded Lie algebra over $\bold Q$ defined 
combinatorially
was studied independently by Y. Ihara (\cite {Ih92} and \cite {Ih99})
and V. G. Drinfel'd (\cite {Dr}).
It can be regarded as a graded Lie algebra version of the Grothendieck-Teichm\"{u}ller group.
An embedding 
\[\varPsi_l:\frak g^l\centerdot \hookrightarrow \frak D\centerdot\otimes_{\bold Q}\bold Q_l,
\]
where $\frak g^l\centerdot$ is the graded Lie algebra over $\bold Q_l$ 
associated to the Galois image of this pro-$l$ representation,
was constructed,
and  it was conjectured to be an isomorphism for all primes $l$,
i.e. $\frak D\centerdot$ may be a common $\bold Q$-structure of $\frak g^l\centerdot$
for all primes $l$ (\cite{Ih90}). \par
In contrast, here we shall construct a Hodge counterpart of this map.
Let $Z_w$ be the $\bold Q$-vector space generated by all multiple zeta values
of indices with weight $w$.
And put $Z_0=\bold Q$.
The {\it MZV algebra} is its formal direct sum;
$Z\centerdot=\underset{w\geqslant 0}{\oplus}Z_w$.
The {\it new-zeta space} is a filtered graded $\bold Q$-vector space which is a
quotient of the ideal $Z_{>2}$;
$NZ\centerdot=Z_{>2}\left/ (Z_{>0})^2 \right.$.
For more details, see \S\S \ref{new-zeta space}.
\begin{mthm}
There is a canonical surjective  $\bold Q$-linear map 
\[
\varPsi_{DR}:\frak D^*\centerdot\twoheadrightarrow NZ\centerdot
\]
as filtered graded $\bold Q$-vector space,
where $\frak D^*\centerdot$ is the filtered graded dual of the stable derivation algebra.
\end{mthm}
There is a standard conjecture on the structure of $\frak D\centerdot$
which is a combination of the conjectures (\S\S \ref{connection})
by Y. Ihara and P. Deligne
concerning the Galois representation on
$\pi_1^l(\bold P^1_{\overline{\bold Q}}-\left\{0,1,\infty\right\})$.
\begin{mconj}
$\frak D\centerdot$ is a free Lie algebra generated by one element in each
degree $m$ $(=3,5,7,\cdots)$.
\end{mconj}

As a corollary of Theorem \ref{mainthm},
we get the following

\begin{mprop}
If we assume Conjecture \ref{Deligne-Ihara}, then we get 
a ``final upper bound''
\footnote{
A.B. Goncharov announced that he has shown the upper-bounding part of  
Dimension Conjecture by the theory of mixed Tate motives in \cite{GonECM}.
Recently another proof was also given by T. Terasoma in \cite{Te}.
}
\[dim_{\bold Q} Z_w \leqslant d_w \]
where $d_w$ is the conjectured value of $dim_{\bold Q} Z_w$
in Dimension Conjecture 
(\cite{Za}, see also \S\S \ref{MZValg}).
\end{mprop}

Modulo Conjecture 2.4.1, Dimension Conjecture is equivalent to

\begin{mmmconj}
$\varPsi_{DR}$ is an isomorphism.
\end{mmmconj}

Particularly this conjecture claims that the dual vector space of 
$NZ\centerdot$
might admit a co-Lie algebra structure deduced from that of  
$\frak D^*\centerdot$.\par

This paper is organized as follows.
The definitions and basic known facts on multiple zeta values will be reviewed
and the new zeta space will be introduced
in \S 1.
\S 2 is devoted to reviewing the stable derivation algebra which appeared
in some works by Y. Ihara on the Galois representation
on $\pi_1^l(\bold P ^1_ {\overline{\bold Q}}- \left\{0,1,\infty\right\})$.
In \S 3, we shall recall the basic properties of the Drinfel'd associator
which are required to prove our main theorem.
We shall state and prove our main results in \S 4.
The connection between the new-zeta space and the stable derivation algebra will be built up there.
In \S 5, we shall compare our story in Hodge side (\S 4) with the story in Galois side (\S 2).

\begin{ack}
I am deeply grateful to
Professor Y. Ihara for introducing me to this area.
This paper cannot be written without his continuous encouragements.
I would also like to thank
H. Tsunogai for kindly showing me  his many computations
and
M. Kaneko for informing me literatures on MZV's.
\end{ack}

\section{The MZV algebra}
This section is devoted to reviewing the known facts and basic conjectures on 
multiple zeta values
and to introducing the new-zeta space.

\subsection{Introduction of MZV}
\begin{defn}
 For each index $\bold{k}=(k_1,k_2,\dotsc ,k_m)$
of positive  integers with
$k_1,\dotsm. k_{m-1}\geqslant 1$, $k_m >1$, 
the corresponding {\sf multiple zeta value} ({\sf MZV} for short) 
$\zeta(\bold{k})$ is, by definition, the real number defined by the convergent series:
\begin{equation*}
\zeta(\bold{k})=
\underset{n_i\in\bold  N}
{
\underset{0<n_1<\dotsm<n_m}{\sum}
}
\frac{1}{n_1^{k_1}\dotsm n_m^{k_m}}
\qquad  .
\end{equation*}
The {\sf weight} of $ \bold{k}:wt(\bold k)$ 
(resp. the {\sf depth} of $ \bold{k}:dp(\bold k)$ )
is defined as $wt(\bold{k)}=k_1+\dotsm +k_m$ (resp. $ dp(\bold{k})=m$).
\end{defn}
Sometimes a multiple zeta value
is called as a Zagier sum, a multiple harmonic series or a poly zeta value.

\begin{eg}
\begin{align*}
&dp=1,wt=k(k\in \bold N_{\geqslant 2}) &:&\zeta(k)\text{ (Riemann zeta value)}  \\
&wt=1,      &:&\text{no MZV's}  \\
&wt=2, dp=1&:&\zeta(2)\\
&wt=3,dp=1&:&\zeta(3)\\
& \qquad \quad    ,dp=2&:&\zeta(1,2)\\
&wt=4,dp=1&:&\zeta(4)\\
& \qquad \quad    ,dp=2&:&\zeta(1,3),\zeta(2,2)\\
& \qquad \quad    ,dp=3&:&\zeta(1,1,2)\\
&wt=w(\geqslant 2),      &:&2^{w-2} \text{ tuples of MZV's}.\\
\end{align*}
\end{eg}

\begin{defn}\label{MZVvec}
For each natural number $w$,
 let $Z_w$ be the $\bold Q$-vector space of $\bold R$ generated by all MZV's of indices with  weight $w$:
$Z_w=\langle\zeta(\bold{k})|wt(\bold{k})=w\rangle_{\bold Q}\subseteq\bold R$,
and put $Z_0=\bold Q$.
For each natural number $m$ and $w$, let
$Z^m_w$ be the subspace of $Z_w$ generated by all MZV's 
of indices with weight $w$ and depth $m$:
$Z^m_w=\langle \zeta(\bold{k})|wt(\bold{k})=w,dp(\bold{k})=m\rangle_{\bold Q}$.\par
$Z_w$ has  the ascending {\it depth filtration}:
$Z_w^{\leqslant m}=\langle\zeta(\bold{k})|wt(\bold{k})=w,dp(\bold{k})\leqslant m\rangle_{\bold Q}$.
\[
\langle\zeta(w)\rangle_{\bold Q}=Z_w^{\leqslant 1}\subseteq Z_w^{\leqslant 2}
\subseteq \dotsm \subseteq Z_w^{\leqslant w-1}=Z_w\subset\bold R
\]\par
On $Z_0$, set 
$Z_0=Z_0^{\leqslant 0}=Z_0^{\leqslant 1}=Z_0^{\leqslant 2}=\dotsm$.
\end{defn}
We should note that $Z_w$ is not graded by depths.
For example, see $\zeta(3)=\zeta(1,2),\zeta(1,3)=\frac{1}{4}\zeta(4)$.

\subsection{The MZV algebra}\label{MZValg}
Let $Z\centerdot$ be the {\it formal direct sum} of $Z_w$ for all $w\geqslant 0$;
$Z\centerdot=\underset{w\geqslant 0}{\oplus}Z_w$. 
Put $Z_{>0}=\underset{w> 0}{\oplus}Z_w$.
They  are $\bold Q$-vector spaces graded by weights.

\begin{lem}
The product of two MZV's, one with weight $a$ and depth $m$,
the other with weight $b$ and depth $l$
belongs to the vector space $Z_{a+b}^{\leqslant l+m}$.
\end{lem}

This follows directly from the definitions.
Considering $Z\centerdot$ as the sub-vector space of the graded algebra
$S\centerdot=\underset{w\geqslant 0}{\oplus}S_w$ where $S_w=\bold R$,
we can give $Z\centerdot$ a structure of algebra.

\begin{propy}
$Z\centerdot$ becomes a filtered graded $\bold Q$-algebra (i.e.
$Z_a^{\leqslant l}\cdot Z_b^{\leqslant m}\subseteq Z_{a+b}^{\leqslant l+m}$)
and $Z_{>0}$ is a homogeneous ideal of $Z\centerdot$.
\end{propy}

We shall call $Z\centerdot$ the {\sf MZV algebra}. 
There is a natural ring homomorphism
$d:Z\centerdot\to\bold R$
which is identity on $Z_w$ for each $w$.
The following conjecture is stated in \cite {Gon97} and \cite{GonECM}.

\begin{dsconj}
The homomorphism $d$ is injective,
i.e. there are no non-trivial 
$\bold Q$-linear relations among different weight MZV's.
\end{dsconj}

This conjecture would imply that all MZV's are transcendental numbers.
But not much is known about these transcendencies.
Among the  known results, we quote here the following;
$\zeta (2k)$ $(k=1,2,\dots)$
is a rational multiple of $\pi^{2k}$ (L. Euler), 
it is transcendental (Lindemann), 
$\zeta (3)$ is irrational (R. Ap\'{e}ry) and 
$\{\zeta (n)\notin \bold Q|n=3,5,7,\dots\}$
is an infinite set, which is a recent result of Tanguy Rivoal \cite{Ri}.
It seems hard to prove the direct sum conjecture.
On the dimension of $Z\centerdot$, the following conjecture appears in \cite{Za}.

\begin{dconj}
$dim_{\bold Q }Z_w$ is  equal to $d_w$, which is given by the Fibonacci-like recurrence $d_w =d_{w-2}+d_{w-3}$, with initial values 
$d_0=1,d_1=0,d_2=1$, 
i.e.
\[
d_w=-\frac{\alpha^{w+2}(\beta-\gamma)+\beta^{w+2}(\gamma-\alpha)+\gamma^{w+2}(\alpha-\beta)}{(\alpha-\beta)(\beta-\gamma)(\gamma-\alpha)}\qquad ,\]
where $\alpha$,$\beta$,$\gamma$ are roots of $x^3-x-1$ 
among which
$\alpha$ is the only real root $\alpha=1.324717\dotsm$,
and $|\beta|=|\gamma|<1$.

\end{dconj}
By the above formula for $d_w$,
approximately 
$d_w \sim \frac{\alpha^{w+2}}
{(\alpha-\beta)(\alpha-\gamma)}=
(0.411496\dotsm)\cdot(1.324717\dotsm)^w$,
so $d_w$ is much smaller than $2^{w-2}$, 
the number of indices of $wt=w$.
So it indicates that there must be many $\bold Q$-linear relations among MZV's
of the same weight.
In fact, several such relations have been 
found by many mathematicians and physicists.
For example,
$\zeta(3)=\zeta(1,2)$ and
$\zeta(4)=4\zeta(1,3)=\frac{4}{3}\zeta(2,2)$
 by L. Euler.
Due to those relations we know 
generators of the MZV algebra in lower weights in Example \ref{MZV} and 
the upper bound  $dim_{\bold Q }Z_w \leqslant d_w$ for $w\leqslant 12$
(for $d_w$, see also below).
But it seems difficult to prove the lower-bounding inequalities
because we need to show the linear independency.

\begin{tabular}[t]{|c||c|c|c|c|c|c|c|c|c|c|c|c|c|}\hline
$w$& 0&1&2&3&4&5&6&7&8&9&10&11&12\\ \hline
$d_w$& 1&0&1&1&1&2&2&3&4&5&7&9&12\\ \hline
$2^{w-2}$ & /& /&1&2&4&8&16&32&64&128&256&512&1024 \\ \hline
\end{tabular}

\begin{eg}\label{MZV}
\begin{align*}
Z_0=&\langle 1 \rangle_{\bold Q} \\
Z_1=&   0  \\
Z_2=&\langle\pi^2\rangle_{\bold Q} \\
Z_3=&\langle\underline{\zeta(3)}\rangle_{\bold Q} \\
Z_4=&\langle\pi^4\rangle_{\bold Q} \\
Z_5=&\langle\pi^2\zeta(3),\underline{\zeta(5)}\rangle_{\bold Q} \\
Z_6=&\langle\pi^6,\zeta(3)^2\rangle_{\bold Q} \\
Z_7=&\langle\pi^4\zeta(3),\pi^2\zeta(5),\underline{\zeta(7)}\rangle_{\bold Q} \\
Z_8=&\langle\pi^8,\pi^2\zeta(3)^2,\zeta(3)\zeta(5),\underline{\zeta(3,5)}\rangle_{\bold Q} \\
Z_9=&\langle\pi^6\zeta(3),\pi^4\zeta(5),\pi^2\zeta(7),\zeta(3)^3,\underline{\zeta(9)}\rangle_{\bold Q} \\
Z_{10}=&\langle\pi^{10},\pi^4\zeta(3)^2,\pi^2\zeta(3)\zeta(5),\pi^2\zeta(3,5),\zeta(3)\zeta(7),\zeta(5)^2,\underline{\zeta(3,7)}\rangle_{\bold Q} \\
\end{align*}
For weights $11$ and $12$, see Example \ref{11&12}.
\end{eg}

\subsection{The new-zeta space}\label{new-zeta space}
It seems now that there exist
no non-trivial $\bold Q$-linear relations among the
above MZV's at each weight in Examples 1.2.3.
So if the dimension conjecture is true, 
they must form bases of the vector spaces in respective weights.
Compare the underlined MZV's with the other ones.
The other ones are old-comers in the sense that
they are written as a 
product of lower weight MZV's.
In contrast, the underlined ones might be new-comers.
Those {\it new-zeta values} are algebraic generators of the MZV algebra $Z\centerdot$.
{\it Where and how many those new-zeta values appear? }
On this question,
it is natural to take the following quotient space of $Z\centerdot$ .
\begin{defn}
The {\sf new-zeta space} is the graded vector space over $\bold Q$.
\[NZ\centerdot=\underset{w\geqslant 1}{\oplus}NZ_w
:=Z\centerdot\left/ \Bigl( (Z_{>0})^2 \oplus  \bold Q\cdot\pi^2  \oplus \bold Q\cdot 1\Bigr) \right.
=Z_{>2}\left/ (Z_{>0})^2 \right. \ ,
\]
where $(Z_{>0})^2$ is the ideal of $Z\centerdot$ generated by products
of the two elements in $Z_{>0}$.
It is equipped with the depth filtration
$\{NZ\centerdot^{\leqslant m}\}_{m\in \bold N}$
which is induced from the depth filtration
$\{Z\centerdot^{\leqslant m}\}_{m\in \bold N}$ of $Z\centerdot$
by the natural surjection
$Z\centerdot\twoheadrightarrow NZ\centerdot$.
\end{defn}
It is also natural to raise the following 

\begin{prob}
Clarify the structure of 
$\Bigl(NZ\centerdot,\{NZ\centerdot^{\leqslant m}\}_{m\in \bold N}\Bigr)$
 as the filtered graded $\bold Q$-vector space.
\end{prob} 

This problem is one of our main motivations in this paper.
We note that the algebra structure of $NZ\centerdot$ is no more interesting,
because its multiplicative map is trivial.

\section{Review of the stable derivation algebra}
In this section, we shall recall the definitions on the
stable derivation algebra which appeared in works by Ihara on 
Galois representation on
$\pi \sb 1(\bold P ^1\sb{\overline{\bold Q}}-\left\{0,1,\infty\right\})$
(\cite{Ih89}$\sim$\cite{Ih99}).

\subsection{The stable derivation algebra}
Let $n\geqslant 4$.
The {\it braid Lie algebra on $n$ strings}
$\frak P^{(n)}_\centerdot=\underset{k\geqslant 1}{\oplus}\frak P^{(n)}_k$
is the graded Lie algebra over $\bold Q$ 
which has the following presentation.

\begin{description}
\item[Generators] \ \ \ $x_{i,j}  \ \ \ (1 \leqslant i,j \leqslant n)$ 
\item[Relations] 
\begin{enumerate}
\renewcommand{\labelenumi}{(\alph{enumi})}
\item  $x_{i,i}=0 \ \ ,x_{i,j}=x_{j,i}  \ \ (1\leqslant i,j \leqslant n)$
\item  $\sum\limits_{k=1}^5 x_{i,k}  =0  \ \ (1\leqslant i\leqslant n)$
\item  $[x_{i,j},x_{k,l}]=0 \ \   \  (\{i,j\} \cap \{k,l\} =\emptyset )$
\end{enumerate}
\item[Grading]  \ \ $ deg \ x_{i,j}=1 \ \ \ \ (1\leqslant i<j \leqslant n)$
\end{description}
It is obtained from the lower central series of the pure braid 
group on $n$ strings in a cetain standard manner (see \cite{Ih91}).
Note that $\frak P^{(4)}_\centerdot$ is freely generated by $x_{1,2}$ and $x_{2,3}$,
and  $\frak P_\centerdot^{(5)}$ is an extension of the free Lie algebra of rank $2$ by that of rank $3$.
As for general $\frak P^{(n)}_\centerdot$,
they are connected with each other  by the standard projections
$p_n:\frak P^{(n+1)}_\centerdot\twoheadrightarrow
\frak P^{(n)}_\centerdot$
defined by
$x_{i,j}\mapsto x_{i,j}$ $(1\leqslant i,j \leqslant n)$ and 
$x_{i,j}\mapsto 0$ (othewise).\par

The derivation $D$ of $\frak P^{(n)}_\centerdot$ is called 
{\it special} if there exists some 
$t_{i,j}\in \frak P^{(n)}_\centerdot$
such that
$D(x_{i,j})=[t_{i,j},x_{i,j}]$ for each 
$i,j$ ($1\leqslant i,j \leqslant n$).
The special derivations of $\frak P^{(n)}_\centerdot$ 
form the graded Lie algebra  $Der^\sharp\frak P^{(n)}_\centerdot$
whose degree $k$ ($k\geqslant 1$) part consists of those $D$'s with 
$t_{i,j}\in \frak P^{(n)}_k$ ($1\leqslant i,j \leqslant n$).
This Lie algebra contains the inner derivations 
$Inn Der^\sharp\frak P^{(n)}_\centerdot$
as homogeneous ideal
and the quotient Lie algebra 
$Out Der^\sharp\frak P^{(n)}_\centerdot
:={Der^\sharp\frak P^{(n)}_\centerdot}/{Inn Der^\sharp\frak P^{(n)}_\centerdot}$
is called the graded Lie algebra of 
{\it special outer derivations} of $\frak P^{(n)}_\centerdot$.
The symmetric group $\frak S_n$ acts on $\frak P^{(n)}_\centerdot$
as automorphisms of graded Lie algebra by 
$x_{i,j}\mapsto x_{\sigma(i),\sigma(j)}$ ($\sigma\in\frak S_n$)
for each $i$, $j$ ($1\leqslant i,j \leqslant n$).
The $\frak S_n$-action
$D\mapsto \sigma\circ D\circ\sigma^{-1}$ for each derivation $D$
of $\frak P^{(n)}_\centerdot$ ($\sigma\in\frak S_n$) induces
an $\frak S_n$-action on 
$Out Der^\sharp\frak P^{(n)}_\centerdot$.
The invariant subalgebra
 $\frak D^{(n)}_\centerdot:=
(Out Der^\sharp\frak P^{(n)}_\centerdot)^{\frak S_n}$
of this action is called the graded Lie algebra of 
{\it symmetric special outer derivations} of 
$\frak P^{(n)}_\centerdot$.
Note that each degree $k$ piece $\frak D^{(n)}_k$ of 
$\frak D^{(n)}_\centerdot=\underset{k\geqslant 1}{\oplus}\frak D^{(n)}_k$
is finite dimensional $\bold Q$-vector space.\par

Since each special derivation $D$ of $\frak P^{(n+1)}_\centerdot$
leaves $Ker p_n$ stable,
we can associate the Lie algebra homomorphism
$\psi_n: \frak D^{(n+1)}_\centerdot\to\frak D^{(n)}_\centerdot$.
On the following sequence of Lie algebra homomorphisms
\[
\dotsm
\overset{\psi_{n+1}}{\to}
\frak D^{(n+1)}_\centerdot
\overset{\psi_n}{\to}
\frak D^{(n)}_\centerdot
\overset{\psi_{n-1}}{\to}
\dotsm\dotsm
\overset{\psi_5}{\to}
\frak D^{(5)}_\centerdot
\overset{\psi_4}{\to}
\frak D^{(4)}_\centerdot ,
\]
Y. Ihara proved that $\psi_n$ is injective for $n\geqslant 4$
in \cite{Ih91}, and then in \cite{Ih92},
that $\psi_n$ is bijective for $n\geqslant 5$.
Therefore, $\psi_n$ induces 
\[
\dotsm
\simeq
\frak D^{(n+1)}_\centerdot
\simeq
\frak D^{(n)}_\centerdot
\simeq
\dotsm\dotsm
\simeq
\frak D^{(5)}_\centerdot
\hookrightarrow
\frak D^{(4)}_\centerdot \ \ .
\]
After the stability property of this tower for $n\geqslant 5$,
he called $\frak D^{(5)}_\centerdot$ the 
{\it stable derivation algebra}.
The sense of considering such a tower of Lie algebras is explained in 
\cite{Ih90} in connection with the action of the absolute Galois group
$Gal({\overline{\bold Q}}/\bold Q)$
on the pro-$l$ fundamental group of
$\bold P^1_{\overline{\bold Q}}-\left\{0,1,\infty\right\}$
and partly  it concerns the philosophy of 
``Teichm\"{u}ller-Lego'' in \cite{Gr}.

\subsection{The normalization form of the 
stable derivation algebra}\label{defSDA}
Let
$\Bbb L\centerdot =\underset{w\geqslant 1}{\oplus}\Bbb L_w
=\bold Q  x\oplus\bold Q y\oplus\bold Q [x,y]\oplus\bold Q [x,[x,y]]\oplus
\bold Q [y,[x,y]]\oplus\cdots$
be the free graded Lie algebra over $\bold Q$ on two variables
$x$ and $y$ with grading $deg \  x=deg \  y=1$.
This Lie algebra $\Bbb L\centerdot$ has
the descending filtration 
$\{\frak F^m\Bbb L\centerdot\}_{m\in\bold N}:
\Bbb L\centerdot=\frak F^1\Bbb L\centerdot \supseteq
\frak F^2\Bbb L\centerdot \supseteq
\frak F^3\Bbb L\centerdot \supseteq
\frak F^4\Bbb L\centerdot \supseteq\dotsm\dotsm$,
which is defined inductively by
$\frak F^1\Bbb L\centerdot=\bold Q y\oplus\bold Q [x,y]\oplus\bold Q [x,[x,y]]\oplus\bold Q [y,[x,y]]\oplus\cdots$ and $\frak F^m\Bbb L\centerdot=[\frak F^1\Bbb L\centerdot,\frak F^{m-1}\Bbb L\centerdot]$.\par

It was shown in \cite{Ih92} that the stable derivation algebra 
$\frak D^{(5)}_\centerdot$  can be identified with  the graded
Lie subalgebra $\frak D_\centerdot$ of $Der \Bbb L\centerdot$
which has the following presentation: \ \ \ 
$\frak D \centerdot=\underset{w\geqslant 1}{\oplus}\frak D_w$ ,
\[
\text{ where  } \ \
\frak D_w=\{D_f\in Der\Bbb L\centerdot| \ f\in\Bbb  L_w
\text {   satisfies {\rm(0)} $\sim$ {\rm(iii)} below.}\}
\]
\[
\begin{cases}
{\rm(0)} \ \   f\in [\Bbb L\centerdot,\Bbb L\centerdot] \ \ (=\overset{\infty}{\underset{k=2}{\oplus}}\Bbb L_k) \\
{\rm(i)} \ \   f(x,y)+f(y,x)=0   \\
{\rm(ii)} \ \  f(x,y)+f(y,z)+f(z,x)=0  \ \ \ \ \ \ \text{for} \ \ x+y+z=0  \\
{\rm(iii)}  \  \sum\limits_{i\in \Bbb Z  /5} f(x_{i,i+1} ,x_{i+1,i+2})=0 \ \ \text{in}  \ \ \frak P^{(5)}_\centerdot \ \ \ .
\end{cases}
\]

Here, for any Lie algebra $H$ and $\alpha,\beta \in H$, $f(\alpha,\beta)$ denotes the image of $f\in\Bbb L\centerdot$ by the homomorphism
$\Bbb L\centerdot \to H$ defined by
$x\mapsto\alpha$, $y\mapsto\beta$ 
and, for $f$ in $\Bbb L\centerdot$, 
$D_f:\Bbb L\centerdot\to \Bbb L\centerdot$ is 
a special derivation defined by $D_f(x)=0$ and $D_f(y)=[y,f]$.
It can be checked easily that $[D_f,D_g]=D_h$ with
$h=[f,g]+D_f(g)-D_g(f)$.
Note that $D\in\frak D\centerdot$  determines uniquely
$f\in [\Bbb L\centerdot,\Bbb L\centerdot]$ such that $D=D_f$.
We also remark that in fact (iii) implies (i).\par

The Lie algebra $\frak D_\centerdot$ was also studied independently
by V. G. Drinfel'd (\cite{Dr}) as 
the graded $\bold Q$-Lie algebra version 
of the {\it Grothendieck-Teichm\"{u}ller group}.
In \cite{Dr}, $\frak D_\centerdot$ and  $\frak P^{(5)}_\centerdot$
 appeared as
$\oplus \frak g\frak r\frak t^n_1(\Bbb Q)$ and $\frak a_4^{\Bbb Q}$.

\begin{rem}
In \cite{Ih92} and \cite{Ih99}, (ii) is replaced by
\begin{equation}\tag*{(ii)$'$}
[y,f(x,y)]+[z,f(x,z)]=0 \ \ \ \ \ \ \text{for} \ \ x+y+z=0.  
\end{equation}
But (i), (ii), (iii) and (i), (ii)$'$, (iii) are equivalent
(cf. \cite{Dr} and \cite{Ih92}).
\end{rem}
From now on, we will identify the stable derivation algebra 
$\frak D^{(5)}_\centerdot$ with
$\frak D_\centerdot$.

\subsection{Weights and depths}
Against the standard terminologies in graded Lie algebras, degrees of 
the stable derivation algebra are called {\it weights}  in
connection with  weights of modular forms in \cite{Ih99}.
The {\it depth filtration} of $\frak D\centerdot$
which is introduced in \cite{Ih99} is
the descending filtration 
$\{\frak F^m \frak D\centerdot\}_{m\in \bold N}$:
\[
\frak D\centerdot=\frak F^1 \frak D\centerdot\supseteq
\frak F^2 \frak D\centerdot\supseteq
\frak F^3 \frak D\centerdot\supseteq\dotsm
\]
\[
\text{ where  }  \ \ \ \
\frak F^m \frak D\centerdot=\{D_f|  \ f\in\frak F^m \Bbb L\centerdot
\text{ satisfies (0) $\sim$ (iii)}\}.
\]
It provides $\frak D\centerdot$ with
a structure of filtered graded Lie algebra.
In \S 4, we shall see that these terminologies, weight and depth,
correspond to those of MZV's.
On the dimension of the depth filtration
of the stable derivation algebra at each weight $w$,
the following results are known.

\begin{prop}[\cite{Ih99}]\label{Takao}
\begin{itemize}
\renewcommand{\labelenumi}{(\roman{enumi})}
\item $\frak F^m \frak D_w=\frak F^{m+1} \frak D_w$  if  
  $w\not\equiv m (mod 2)$.
\item  $\frak F^m \frak D_w=\{0\}$   if $m>\frac{w}{2}$.
\item  $dim_{\bold Q} \ (\frak F^1\frak D_w\left/\frak F^2\frak D_w)\right.=
\begin{cases}
1 & w=3,5,7,9,\dots    \\
0 & w:\text{otherwise   .}
\end{cases} 
$
\item $dim_{\bold Q} \ (\frak F^2\frak D_w\left/\frak F^3\frak D_w)\right.=
\begin{cases}
\  \ \ 0                & w:\text{odd}  \  \\  
 \ [\frac{w-2}{6}]  & w:\text{even}   \ \ \ \text{{\rm(Ihara-Takao)}}.
\end{cases} $
\end{itemize}
\end{prop}
\begin{eg}\label{Matsumoto}
The following computation table of 
the stable derivation algebra up to weight $12$ 
is due to Y. Ihara, M. Matsumoto and H. Tsunogai.
\begin{align*}
&\frak D_1=0,   \\
&\frak D_2=0,   \\
&\frak D_3=\langle D_{f_3}\rangle_\bold Q,   \\
&\frak D_4=0,   \\
&\frak D_5=\langle D_{f_5}\rangle_\bold Q,   \\
&\frak D_6=0,   \\
&\frak D_7=\langle D_{f_7}\rangle_\bold Q,   \\
&\frak D_8=\langle  \ [D_{f_3},D_{f_5}] \  \rangle_\bold Q,   \\
&\frak D_9=\langle D_{f_9}\rangle_\bold Q,   \\
&\frak D_{10}=\langle  \ [D_{f_3},D_{f_7}] \  \rangle_\bold Q,   \\
&\frak D_{11}=\langle  \ D_{f_{11}} ,[D_{f_3}, [D_{f_3},D_{f_5}]] \  \rangle_\bold Q,   \\
&\frak D_{12}=  \langle  \ [D_{f_3},D_{f_9}] ,[D_{f_5},D_{f_7}] \  \rangle_\bold Q.   \\
\end{align*}
Here, $f_w$ is a certain element of $\Bbb L_w$.
For example, $f_3$ and $f_5$ are as given below.

\begin{align*}
f_3 = & [x,[x,y]]+[y,[x,y]],    \\
f_5= & 2[x,[x,[x,[x,y]]]]+4[y,[x,[x,[x,y]]]]+4[y,[y,[x,[x,y]]]] \\
& +2[y,[y,[y,[x,y]]]]+[[x,y],[x,[x,y]]]+3[[x,y],[y,[x,y]]].
\end{align*}

\begin{tabular}[t]{|c||c|c|c|c|c|c|c|c|c|c|c|c|}\hline
$w$ & 1&2&3&4&5&6&7&8&9&10&11&12\\ \hline
$dim_{\bold Q} \ \frak D_w(=\frak F^1\frak D_w)$& 0&0&1&0&1&0&1&1&1&1&2&2\\ 
\hline\hline
$dim_{\bold Q} \ \frak F^2\frak D_w$& 0&0&0&0&0&0&0&1&0&1&1&2\\ \hline
$dim_{\bold Q} \ \frak F^3\frak D_w$& 0&0&0&0&0&0&0&0&0&0&1&1\\ \hline
$dim_{\bold Q} \ \frak F^4\frak D_w$& 0&0&0&0&0&0&0&0&0&0&0&1\\ \hline
$dim_{\bold Q} \ \frak F^5\frak D_w$& 0&0&0&0&0&0&0&0&0&0&0&0\\ \hline
\end{tabular}

A remarkable phenomena are pointed out in \cite{Ih89}
(see also \cite{Ih99}).
\[
2[D_{f_3},D_{f_9}]-27[D_{f_5},D_{f_7}]\in\frak F^3\frak D.
\]
\end{eg}

\subsection{A connection between the  $l$-adic Galois image Lie algebra and the stable derivation algebra}\label{connection}
Let $l$ be a prime.
Consider the outer action
\[
\varphi_l:Gal({\overline{\bold Q}}/\bold Q)\to
Out \ \pi_1^{(l)}(\bold P^1_{\overline{\bold Q}}-\left\{0,1,\infty\right\})
\]
of the absolute Galois group $Gal({\overline{\bold Q}}/\bold Q)$
on the pro-$l$ fundamental group of
$\bold P^1_{\overline{\bold Q}}-\left\{0,1,\infty\right\}$.
The {\it $l$-adic Galois image Lie algebra} 
$\frak g^l_\centerdot=\underset{w\geqslant 1}{\oplus}\frak g^{l}_w$
is the graded Lie algebra over $\bold Q_l$
which is constructed from the image of $\varphi_l$
in a certain standard way and it is shown that $\frak g^l_\centerdot$
can be embedded into the $l$-adic stable derivation algebra
\[
\varPsi_l:\frak g^l_\centerdot\hookrightarrow\
\frak D_\centerdot\otimes_\bold Q \bold Q_l
\]
for all primes $l$ (for more details, see \cite{Ih90}).
On the map $\varPsi_l$ and
these Lie algebras, $\frak g^l\centerdot$ and $\frak D\centerdot$,
the following conjectures  are stated in \cite{Ih90} and \cite{Ih99}
and partly in \cite{De}.
\begin{conj}\label{Deligne-Ihara}
\begin{enumerate}
\renewcommand{\labelenumi}{(\arabic{enumi})}
\item $\varPsi_l:\frak g^l_\centerdot\hookrightarrow\frak D_\centerdot\otimes_\bold Q \bold Q_l$ is an isomorphism for every prime $l$.
\item These Lie algebras are free.
\item On generators of these  Lie algebras,
we can take one element in each odd degree $m \ (m=3,5,7,\cdots)$.
\end{enumerate}
\end{conj}

The validity of (1) implies that the $l$-adic Galois image Lie algebra 
$\frak g^l\centerdot$ has a common $\bold Q$-structure
for all primes $l$ and 
the stable derivation algebra
$\frak D\centerdot$ which is defined entirely
combinatorially becomes this common $\bold Q$-structure.
Conjecture (3) for $\frak g^l\centerdot$  was proved by R. Hain and M. Matsumoto 
in \cite{HaMa}.
On the lower weight up to $12$ for $\frak D\centerdot$ and
$17$ for $\frak g^l\centerdot$,
H. Tsunogai verified (1), (2) and (3)  by computations (\cite{Tsu}).

\begin{tabular}[t]{|c||c|c|c|c|c|c|c|c|c|c|c|c|c|c|c|c|c|}\hline
$w$ & 1&2&3&4&5&6&7&8&9&10&11&12&13&14&15&16&17\\ \hline
$dim_{\bold Q} \ \frak D_w$ & 0&0&1&0&1&0&1&1&1&1&2&2&/&/&/&/&/ \\ \hline
$dim_{\bold Q} \ \frak g^l_w$ & 0&0&1&0&1&0&1&1&1&1&2&2&3&3&4&4\text{ or }5
&7 \\ \hline
\end{tabular}

\section{Review of the Drinfel'd associator}
We shall review the basic properties of the Drinfel'd associator which is required to prove our main theorem in 
the next section. Most of them can be found in \cite{Dr} and \cite{Kas}.

\subsection{The KZ equation and the Drinfel'd associator}
Let $A^{\land}_{\Bbb C}=\Bbb C\langle\langle A,B\rangle\rangle$
be the non-commutative formal power series ring with complex number 
coefficients
generated by two elements $A$ and $B$.
Consider the {\it Knizhnik-Zamolodchikov equation} 
({\it KZ equation} for short)
\begin{equation}
\tag*{(KZ)}
\frac{\partial G}{\partial u}(u)=
(\frac{A}{u}+\frac{B}{u-1})\cdot G(u) \ \ \ \ \ \ \ \ \ \ ,
\end{equation}
where $G(u)$ is an analytic function in complex variable $u$
with values in  $A^{\land}_{\Bbb C}$
where \lq analytic' means  
each of whose coefficient is analytic.
The equation (KZ) has singularities only at $0,1$ and $\infty$.
Let $\Bbb C'$ be the complement of the union of the real 
half-lines $(-\infty,0]$ and $[1,+\infty)$ in the complex plane.
This is a simply-connected domain.
The equation (KZ) has a unique analytic solution on $\Bbb C'$ having a specified value 
at any given points on  $\Bbb C'$.
Moreover, at the singular points $0$ and $1$,
there exist unique solutions $G_0(u)$ and $G_1(u)$ of (KZ) such that
\[
G_0(u)\approx u^A  \  (u\to 0), \ \   \ \ \ \  G_1(u)\approx (1-u)^B  \  (u\to 1),  
\]
where $\approx$ means that 
$G_0(u)\cdot u^{-A}$ (resp. $G_1(u)\cdot(1-u)^{-B}$) 
has an analytic  continuation in a neighborhood of $0$ (resp. $1$) 
with value $1$ at $0$ (resp. $1$).
Here,
$u^A:=exp(A log u)
:=1+\frac{A log u}{1!}+\frac{(A log u)^2}{2!}+\frac{(A log u)^3}{3!}+\dotsb$
and
$log u := \int_1^u\frac{dt}{t}$ in $\Bbb C'$.
In the same way, $(1-u)^B$ is well-defined on $\Bbb C'$.
One can calculate the lower degree parts of $G_0(u)$ as follows.

\begin{align*}
G_0(z)&= 1+(log z)A+{log (1-z)}B  +\frac{(log z)^2}{2}A^2-Li_2(z)AB\\ 
& +\left\{Li_2(z)+(log z) log(1-z)\right\}BA 
+\frac{\{log (1-z)\}^2}{2}B^2
+\frac{(log z)^3}{6}A^3\\
&-Li_3(z)A^2B
  +\left\{2Li_3(z)+(log z) Li_2(z)\right\}ABA 
+Li_{1,2}(z)AB^2 \\
& -\left[Li_3(z)-(log z) Li_2(z)-\frac{(log z)^2log(1-z)}{2}\right]BA^2 
+Li_{2,1}(z)BAB \\
&-\left[Li_{1,2}(z)+Li_{2,1}(z)-\frac{log z\{log(1-z)\}^2}{2}\right]B^2A \\
&+\frac{\{log (1-z)\}^3}{6}B^3+\cdots \\
\end{align*}
where
\[
Li_{k_1,\dotsc ,k_m}(z):=
\underset{n_i\in\bold  N}
{
\underset{0<n_1<\dotsm<n_m}{\sum}
}
\frac{z^{n_m}}{n_1^{k_1}\dotsm n_m^{k_m}} \  \ .
\]

Since $G_0(u)$ and $G_1(u)$ are both non-zero unique solutions of (KZ) 
with the specified asymptotic behaviors, 
they must coincide with each other up to multiplication
by an invertible element of $A^{\land}_{\Bbb C}$.

\begin{defn}
The {\sf Drinfel'd associator}
\footnote{
To be precise, Drinfel'd defined $\varphi_{KZ}(A,B)$ instead of $\Phi_{KZ}(A,B)$ in \cite{Dr}, 
where $\varphi_{KZ}(A,B)=\Phi_{KZ}(\frac{1}{2\pi i}A,\frac{1}{2\pi i}B)$.}
is the  element $\Phi_{KZ}(A,B)$
of $A^{\land}_{\Bbb C}$ which is defined  by
\[
\ \ \ \ G_0(u)=G_1(u)\cdot\Phi_{KZ}(A,B) \  \ .
\]
\end{defn}

By considering the image in $(A^{\land}_{\Bbb C})^{ab}$, 
the abelianization of
$A^{\land}_{\Bbb C}$, we easily find that $\Phi_{KZ}(A,B)\equiv 1$ in
$(A^{\land}_{\Bbb C})^{ab}$.

\subsection{Explicit formulae}\label{Property I}
We will discuss on each coefficient 
of the Drinfel'd associator $\Phi_{KZ}(A,B)$.\par
Let $\omega_1,\omega_2,\dotsc,\omega_n$ $(n\geqslant 1)$
be differential $1$-forms on $\Bbb C'$.
An iterated integral 
$\int_0^1\omega_n\circ\omega_{n-1}\circ\dotsi\circ\omega_1$
is defined inductively as
$\int_0^1\omega_n(t_n)\int_0^{t_n}\omega_{n-1}
\circ\dotsi\circ\omega_1$.
It is known that MZV's can be written by iterated integrals as follows.
\begin{align*}
\zeta(k_1,k_2,\dotsc,k_m)
= \int_0^1 
\underbrace{\frac{du}{u}\circ\dotsm\circ\frac{du}{u}\circ
 \frac{du}{1-u}}_{k_m}& \circ 
\frac{du}{u}\circ\dotsm\dotsm 
\circ\frac{du}{1-u} \\& 
\circ
\underbrace{\frac{du}{u}\circ\dotsm\circ\frac{du}{u}\circ\frac{du}{1-u}}_{k_1}
\end{align*}
This expression is due to Kontsevich and Drinfel'd.
It can be verified by direct calculations
(see, for example, \cite{Gon97} and \cite{GonECM}).

Let 
\[
\Bbb A\centerdot=\underset{w\geqslant 0}{\oplus}\Bbb A_w=\bold Q\langle A,B\rangle (\subset A^{\land}_{\Bbb C})
\]
be the non-commutative graded polynomial ring 
over $\bold Q$
with two variables $A$ and $B$ with
$deg A= deg B =1$.
Here $\Bbb A_w$ is the homogeneous degree $w$ part of $\Bbb A\centerdot$.
Put 
\[
M=A\cdot\Bbb A\centerdot\cdot B=\{A\cdot F\cdot B| F\in\Bbb A\centerdot\}
\]
which is the $\bold Q$-linear subspace of 
$\Bbb A\centerdot$.
Define the $\bold Q$-linear map
$Z:M\to \Bbb C$ which is determined by
\begin{align*}
Z(A^{p_1}B^{q_1} & A^{p_2}B^{q_2}\dotsm A^{p_k}B^{q_k})
:= \int_0^1 
\underbrace{\frac{du}{u}\circ\dotsm\circ\frac{du}{u}}_{p_1}
\circ
\underbrace{\frac{du}{1-u}\circ\dotsm\circ\frac{du}{1-u}}_{q_1} \\& 
\circ
\frac{du}{u}\circ\dotsm\dotsm\dotsm
\circ\frac{du}{1-u} 
\circ
\underbrace{\frac{du}{u}\circ\dotsm\circ\frac{du}{u}}_{p_k}
\circ
\underbrace{\frac{du}{1-u}\circ\dotsm\circ\frac{du}{1-u}}_{q_k}  \\ & \qquad
=\zeta(
\underbrace{1,\ldots 1}_{q_k-1}
,p_k+1,
\underbrace{1,\ldots 1}_{q_{k-1}-1}
,p_{k-1}+1,\ldots\ldots ,1,p_1+1)
\end{align*}
for $p_i,q_i\geqslant 1 (1\leqslant i \leqslant k)$.
It is the iterated integral from $0$ to $1$ obtained by replacing 
$A$ by $\frac{du}{u}$ and $B$ by $\frac{du}{1-u}$.
\begin{defn}
The {\sf word} is an element of $\Bbb A\centerdot$
which is monic and monomial.
But exceptionally we shall not call $1$ a word.
For each word $W$, the {\sf weight} and {\sf depth} of $W$ are as follows.\par
$wt (W):=$ `the sum of exponents of $A$ and $B$ in $W$' \par
$dp (W):=$ `the sum of exponent  \ of\ \ \ \ \ \ \ \ \ \ $B$ in $W$'
\end{defn}
For example, $A^3BAB$ is a word with $wt (A^3BAB)=6$ and $dp (A^3BAB)=2$.
Note that, by definition, $Z(W)\in Z^m_w$ (see Definition \ref{MZVvec})
 for any word $W$ with 
$wt (W)=w$ and $dp( W)=m$.
We can expand uniquely as
$\Phi_{KZ}(A,B)=1+\sum\limits_{W:\text{words}} I(W)W $,
where each $I(W)$ is a complex number.
The following product structure is required to calculate each $I(W)$
explicitly.

\begin{defn}
The {\sf shuffle product} `$\circ$'
is the $\bold Q$-bilinear map
$\circ:\Bbb A\centerdot\times\Bbb A\centerdot\to\Bbb A\centerdot$,
which is defined by
\begin{description}
\item[\underline{S1}]
$W\circ 1=1\circ W=W$
\item[\underline{S2}]
$UW\circ VW'=U(W\circ VW')+V(UW\circ W') \text{ \ \ where }
U,V\in\{A,B\},$
\end{description}
for all words  $W$, $W'$ in $\Bbb A\centerdot$.
\end{defn}
For example,
$AB\circ A=2A^2B+ABA,AB\circ AB=2ABAB+4A^2B^2$.
Note that 
$Z(W)\cdot Z(W')=Z(W\circ W')$  \ \  if $W,W' \in M$.
It is a basic property of iterated integral.
There is a natural surjection from
$\Bbb A\centerdot$  to
$\Bbb A\centerdot\left/ \Bigl( B \Bbb A\centerdot+\Bbb A\centerdot A\Bigr) \right.$.
By identifying the latter space with $\bold Q\cdot 1+M (=\bold Q\cdot 1+A\cdot\Bbb A\centerdot\cdot B)$
we obtain the $\bold Q$-linear map
$f:\Bbb A\centerdot\hookrightarrow \Bbb A\centerdot\left/ \Bigl( B \Bbb A\centerdot+\Bbb A\centerdot A\Bigr) \right. \overset{\sim}{\to}\bold Q\cdot 1+M \hookrightarrow \Bbb A\centerdot$.

\begin{prop}[(Explicit Formulae)]\label{Explicit Formula}
Each coefficient of 
$\Phi_{KZ}(A,B)=1+\sum\limits_{W:\text{words}} I(W)W $
can be expressed as follows.
\footnote{
Another explicit formula of the Drinfel'd associator was obtained in \cite{LM},
but there seems to be an error on the signature in Theorem A.9.
Their formula is inconvenient for our present purpose
because they expressed each coefficient of MZV in terms of words
instead of expressing each coefficient of word in terms of MZV's.
}
\begin{enumerate}
\renewcommand{\labelenumi}{(\alph{enumi})}
\item  When $W$ is in $M$, \ \ \ $I(W)=(-1)^{dp (W)}Z(W)$.
\item When $W$ is written as $B^rVA^s (r,s\geqslant 0,V\in M)$,
\[
I(W)=(-1)^{dp (W)}\sum_{0\leqslant a\leqslant r,0\leqslant b\leqslant s}(-1)^{a+b}
Z\Bigl(f(B^a\circ B^{r-a}VA^{s-b}\circ A^b)\Bigr).
\]
\item When $W$ is written as $B^rA^s (r,s\geqslant 0)$,
\[
I(W)=(-1)^{dp (W)}\sum_{0\leqslant a\leqslant r ,0\leqslant  b\leqslant s}(-1)^{a+b}
Z\Bigl(f(B^a\circ B^{r-a}A^{s-b}\circ A^b)\Bigr).
\]
\end{enumerate}
\end{prop}
\begin{pf}
Applying the method (A.15) in \cite{LM},
we get these expressions.
\qed
\end{pf}

This proposition implies 

\begin{propyI}
Each coefficient of $\Phi_{KZ}(A,B)$ can be written by MZV's.
More precisely, for each word $W$ with $wt (W)=w$ and $dp (W)=m$,
$I(W)$ is in $Z^m_w$.
\end{propyI}

The terms of low degrees of $\Phi_{KZ}(A,B)$ are as follows.
\begin{align*}
\Phi_{KZ}& (A,B)= 1-\zeta(2)[A,B]-\zeta(3)[A,[A,B]]+\zeta(1,2)[[A,B],B]   \\
&+\zeta(4)[A,[A,[A,B]]]+\zeta(1,3)[A,[[A,B],B]]-\zeta(1,1,2)[[[A,B],B],B] \\
&+\frac{1}{2}\zeta(2)^2[A,B]^2+\dotsm\dotsm.
\end{align*}

\subsection{Relations}\label{Property II}
Let $\Bbb L^{\land}_{\Bbb C}$ be the completion by degree of the free Lie
algebra $\Bbb L\centerdot\otimes\Bbb C$ over $\Bbb C$
(for $\Bbb L\centerdot$, see \S\S \ref{defSDA}).
The Lie algebra $\Bbb L^{\land}_{\Bbb C}$ can be naturally identified with a subspace of 
$\Bbb A^{\land}_{\Bbb C}$
by replacing $[A,B]$ by $AB-BA$.
The following property was a key to prove one of the main theorems in \cite{Dr}.
\begin{propyII}
The Drinfel'd associator satisfies
\[
\begin{cases}
\text{(0)  } log \,\Phi_{KZ}(A,B)
:=\{\sum I(W)W \}
-\frac{1}{2}\{\sum I(W)W \}^2
+\frac{1}{3}\{\sum I(W)W \}^3 \\
\qquad\quad -\frac{1}{4}\{\sum I(W)W \}^4
+\dotsm  
\in [\Bbb L^{\land}_{\Bbb C},\Bbb L^{\land}_{\Bbb C}] \ \ 
\Bigl(={\underset{a\geqslant 2}{\widehat{\oplus}}}(\Bbb L_a\otimes_{\bold Q}\Bbb C)\Bigr) \\
\text{(I)  } \Phi_{KZ}(A,B)\Phi_{KZ}(B,A)=1  \\ 
\text{(II)  }
e^{\pi iA} \Phi_{KZ}(C,A)e^{\pi iC} \Phi_{KZ}(B,C)e^{\pi iB}\Phi_{KZ}(A,B)=1 \\
 \ \ \ \   \ \ \ \   \ \ \ \  \ \ \ \   \ \ \ \   \ \ \ \   \ \ \ \   \ \ \ \  
 \text{                          for  }\ \ \  A+B+C=0 \\
\text{(III)   }
\Phi_{KZ}(x_{1,2} ,x_{2,3})
\Phi_{KZ}(x_{3,4} ,x_{4,5})
\Phi_{KZ}(x_{5,1} ,x_{1,2})  \\
 \ \ \ \ \ \ \  \ \  \ \ 
  \Phi_{KZ}(x_{2,3} ,x_{3,4})
\Phi_{KZ}(x_{4,5} ,x_{5,1})
=1      \ \  \ 
\text{in} \ \ \ \widehat{U\frak P_\centerdot^{(5)}\otimes \Bbb C} \ \ .  \\
\end{cases}
\]
\end{propyII}

Here, $\widehat{U\frak P_\centerdot^{(5)}\otimes \Bbb C}$ 
stands  for the completion by degree
of the universal enveloping algebra of $\frak P_\centerdot^{(5)}$ tensored by $\Bbb C$.
In fact (III) implies (I).

\begin{prop}
Relation $\mathrm{(III)}$ implies $\mathrm{(I)}$.  
\end{prop}
\begin{pf}
Observe that we have a basic projection of completed non-commutative algebras
$p:\widehat{U\frak P_\centerdot^{(5)}\otimes\Bbb C}\twoheadrightarrow\Bbb A^{\land}_{\Bbb C}$
which sends $x_{1,2}$,$x_{2,3}$,$x_{3,4}$,$x_{4,5}$ and $x_{5,1}$
to $A$,$B$,$A$,$0$ and $0$ respectively.
Since
\begin{align*}
p\Bigl( 
\Phi_{KZ}&(x_{1,2} ,x_{2,3})
\Phi_{KZ}(x_{3,4} ,x_{4,5})
\Phi_{KZ}(x_{5,1} ,x_{1,2})
\Phi_{KZ}(x_{2,3} ,x_{3,4})
\Phi_{KZ}(x_{4,5} ,x_{5,1})
\Bigr) \\
&=\Phi_{KZ}(A,B)\Phi_{KZ}(B,A),
\end{align*}
(III) implies (I), as desired. 
\qed
\end{pf}

\begin{rem}
Our above formulae are not exactly the same as those below given in
Drinfel'd's original paper \cite{Dr},
however they are equivalent.

\begin{itemize}
\item 
$log \,\varphi_{KZ}(A,B)\in [\Bbb L^{\land}_{\Bbb C},\Bbb L^{\land}_{\Bbb C}]
 \ \ $
\item 
$\varphi_{KZ}(B,A)=\varphi_{KZ}(A,B)^{-1}  \qquad\qquad\qquad\qquad\qquad\qquad
\quad \ \ \ \ (2.12)$
\item
$ e^{\frac{A}{2}} \varphi_{KZ}(C,A)e^{\frac{C}{2}}\varphi_{KZ}(B,C)e^{\frac{B}{2}}\varphi_{KZ}(A,B)=1$ \\ 
\text{\qquad\qquad\qquad\qquad\qquad\qquad for  }\ $  A+B+C=0 
\qquad\qquad \  \ \ \ \ (5.3)$
\item 
$\varphi_{KZ}(x_{1,2},x_{2,3}+x_{2,4})\varphi_{KZ}(x_{1,3}+x_{2,3},x_{3,4})
=\varphi_{KZ}(x_{2,3},x_{3,4})\cdot $ \par $
 \ \ \ \ \ \  \varphi_{KZ}(x_{1,2}+x_{1,3},x_{2,4}+x_{3,4})
\varphi_{KZ}(x_{1,2},x_{2,3}) \ \ 
\text{   in} \ \ \ \widehat{U\frak P_\centerdot^{(5)}\otimes \Bbb C} \ \ (2.13) $
\end{itemize}
\end{rem}
The proof of his formulae relies on some asymptotic behaviors of 
certain solutions of the KZ equations system on some specified zones.
Property I and formulae (0) $\sim$ (III) suggest that 
we can get many algebraic relations among MZV's.
In the next section, we shall see that his formulae (equivalently ours)
play an essential role to prove our main theorem.

\section{Main results }
The purpose of this section is to show the close relationship between 
$NZ\centerdot$ and $\frak D\centerdot$.

\subsection{The new-zeta quotient of the Drinfel'd associator}
By Property I (\S\S \ref{Property I}), 
the Drinfel'd associator
$\Phi_{KZ}(A,B)=1+\sum\limits_{W:{\text words}} I(W)W $ 
can be regarded as an element of 
${\underset{w\geqslant 0}{\widehat{\oplus}}}
(Z_w\otimes_{\bold Q}\Bbb A_w)$ .

\begin{defn}
The {\sf new-zeta quotient of the Drinfel'd associator}
$\overline{\Phi_{KZ}(A,B)}=\sum\limits_{W:{\text words}}\overline{I(W)}W $
is an element of the $\bold Q$-linear vector space
${\underset{w\geqslant 0}{\widehat{\oplus}}}
(NZ_w\otimes_{\bold Q}\Bbb A_w)$
obtained from $\Phi_{KZ}(A,B)$ by replacing each $I(W)$ by
$\overline{ I(W)}$ which is the image of $I(W)$ by
the natural surjection
$Z_\centerdot\twoheadrightarrow NZ_\centerdot 
(=Z\centerdot\left/ \Bigl((Z_{>0})^2 \oplus Z_2 \oplus Z_0\Bigr)\right. )$.
\end{defn}

Embed the Lie algebra $\Bbb L\centerdot$ into $\Bbb A\centerdot$
by the map $[A,B]\mapsto AB-BA$ and identify $\frak D_w$ with
\[
\left\{
f\in\Bbb A_w 
\left\vert
\begin{array}{cl}
\text{(0)}& \ \   f\in[\Bbb L\centerdot,\Bbb L\centerdot] \ 
(=\overset{\infty}{\underset{k=2}{\oplus}}\Bbb L_k) \\
\text{(i)}& \ \   f(x,y)+f(y,x)=0   \\
\text{(ii)}& \ \  f(x,y)+f(y,z)+f(z,x)=0  \ \ \text{for}\  x+y+z=0  \\
\text{(iii)} & \  \sum\limits_{i\in \Bbb Z  /5} f(x_{i,i+1} ,x_{i+1,i+2})=0 \ \ \text{in}  \ \ \frak P_\centerdot^{(5)}
\end{array}
\right.
\right\}
\]
by sending $D_g\in\frak D_w$ to $g\in\Bbb A_w$ for each $w\geqslant 1$.
Then
we can regard $\frak D\centerdot$ as a subspace of 
the graded vector space $\Bbb A\centerdot$.

\begin{thm}
The new-zeta quotient
of the Drinfel'd associator 
$\overline{\Phi_{KZ}(A,B)}$ lies on the $\bold Q$-linear vector space
${\underset{w\geqslant 2}{\widehat{\oplus}}}
(NZ_w\otimes_{\bold Q}\frak D_w)$.
\end{thm}

\begin{pf}
\begin{align*}
\Phi_{KZ}&=exp \ log\Phi_{KZ}=1+log\Phi_{KZ}+\frac{(log\Phi_{KZ})^2}{2!}
+\frac{(log\Phi_{KZ})^3}{3!}+\cdots \\
&\equiv 1+log\Phi_{KZ} \qquad\qquad\qquad    mod (Z_{>0})^2 .
\end{align*}
This means that $\overline{\Phi_{KZ}}=\sum\limits_{W:{\text words}}\overline{I(W)}W $  lies on 
${\underset{w\geqslant 2}{\widehat{\oplus}}}
(NZ_w\otimes_{\bold Q}\Bbb L_w)$.
From relations (0) $\sim$ (III) in Property II (\S\S \ref{Property II}) of the Drinfel'd associator
$\Phi_{KZ}(A,B)=1+\sum\limits_{W:{\text words}} I(W)W \in
{\underset{w\geqslant 0}{\widehat{\oplus}}}
(Z_w\otimes_{\bold Q}\Bbb A_w)$,
we find that 
$\overline{\Phi_{KZ}(A,B)}=\sum\limits_{W:{\text words}}\overline{I(W)}W \in
{\underset{w\geqslant 2}{\widehat{\oplus}}}
(NZ_w\otimes_{\bold Q}\Bbb L_w)$ satisfies
\begin{align*}
\text{(i)}& \ \   \overline{\Phi_{KZ}(A,B)}+\overline{\Phi_{KZ}(B,A)}=0   \\
\text{(ii)}& \ \  \overline{\Phi_{KZ}(A,B)}+\overline{\Phi_{KZ}(B,C)}+\overline{\Phi_{KZ}(C,A)}=0  \ \ \text{  for  }\  
A+B+C=0  \\
\text{(iii)} & \  \sum\limits_{i\in \Bbb Z  /5}  
\overline{\Phi_{KZ}(x_{i,i+1} ,x_{i+1,i+2})}=0 \quad
\text{ in } \ {\underset{w\geqslant 1}{\widehat{\oplus}}}
(NZ_w\otimes_{\bold Q}\frak P_w^{(5)})
\end{align*}
So we find that $\overline{\Phi_{KZ}}$
can be regarded as an element of
${\underset{w\geqslant 2}{\widehat{\oplus}}}
(NZ_w\otimes_{\bold Q}\frak D_w)$.
\qed
\end{pf}

\subsection{Main theorem}\label{Main Theorem}
\begin{defn}
The {\sf graded dual of the
stable derivation algebra} $\frak D^*_\centerdot=\underset{w\geqslant 1}{\oplus}
\frak D^*_w$
is the graded vector space whose component of degree $w$ is the dual vector
space of $\frak D_w$, i.e.
$\frak D^*_w=(\frak D_w)^*$.
The {\sf depth filtration} 
$\{\frak F^m\frak D^*_\centerdot\}_{m\in \bold N}$
of $\frak D^*_\centerdot$
is the ascending filtration:
\[
\{0\}=\frak F^0\frak D^*_\centerdot\subseteq
\frak F^1\frak D^*_\centerdot\subseteq
\frak F^2\frak D^*_\centerdot\subseteq
\frak F^3\frak D^*_\centerdot\subseteq\dotsm,
\]
where $\frak F^m\frak D^*_w$ is the vector space of linear forms on
$\frak D_w$ whose restriction  to $\frak F^{m+1}\frak D_w$ is zero, i.e.
$\frak F^m\frak D^*_w=\{f^*\in\frak D^*_w \bigl\vert \quad
f^*|_{\frak F^{m+1}\frak D_w}\equiv 0\}$.
\end{defn}
Then 
$\Bigl(\frak D^*_\centerdot,\{\frak F^m\frak D^*_\centerdot\}_{m\in \bold N}\Bigr)$ becomes a filtered graded $\bold Q$-vector space.

\begin{thm}\label{mainthm}
There is a canonical surjective $\bold Q$-linear map 
of filtered graded vector spaces
\[
\varPsi_{DR}:\Bigl(\frak D^*_\centerdot,\{\frak F^m\frak D^*_\centerdot\}_{m\in \bold N}\Bigr)
\twoheadrightarrow 
\Bigl(NZ\centerdot,\{NZ\centerdot^{\leqslant m}\}_{m\in\bold N}\Bigr).
\]
Moreover, it strictly preserves the depth filtration, i.e.
$\varPsi_{DR}(\frak F^m\frak D_w^*)=NZ_w^{\leqslant m}$.
\end{thm}

\begin{pf}
{\it Construction}:
Decompose $\overline{\Phi_{KZ}}$ as
$\overline{\Phi_{KZ}}=\sum\limits_{w=2}^\infty \psi_w$ 
where $\psi_w\in NZ_w\otimes_{\bold Q}\frak D_w$ for all $w$.
There is the natural isomorphism
$NZ_w\otimes_{\bold Q}\frak D_w \cong Hom_{\bold Q}(\frak D_w^*,NZ_w)$,
thus each $\psi_w$ determines a $\bold Q$-linear map
$\Psi_w:\frak D_w^* \to NZ_w$.
Then we define $\varPhi_{DR}$ as follows.
\[
\overset{\infty}{\underset{w=2}{\oplus}}\Psi_w:
\frak D^*_\centerdot\Bigl (=\overset{\infty}{\underset{w=2}{\oplus}}\frak D^*_w\Bigr)
\longrightarrow
NZ_\centerdot\Bigl (=\overset{\infty}{\underset{w=2}{\oplus}}NZ_w\Bigr) \  \ .
\]\par
{\it Surjectivity}:
Let $w\geqslant 2$.
Let $\{W^* \ | \ \ W:\text{words},\ wt(W)=w \}$
be a basis of $\Bbb A^*_w$ (: 
the dual vector space of $\Bbb A_w$) where
\[
W^*(W')=
\begin{cases}
1 & \text{  if   } \ \ W=W', \\
0 & \text{  if   } \ \ W\not= W' ,
\end{cases}
\]
for all words $W' \in \Bbb A_w$.
Denote $W^*|_{\frak D_w}$ be the restriction of the map $W^*$ on $\Bbb A_w$
into $\frak D_w$.
Then it is clear that
$\{W^*|_{\frak D_w} \ | \ W:\text{words}, \ wt(W)=w\}$ 
is a system of generators of the $\bold Q$-vector space $\frak D^*_w$,
i.e.
\[
\frak D^*_w=\left\langle W^*|_{\frak D_w}\Bigl |  \ W:\text{words}, \ wt(W)=w
\right\rangle_{\bold Q}.
\]
Note that 
$\Psi_w(W^*|_{\frak D_w})=\overline{I(W)}$ , thus we find that
\begin{align*}
\Psi_w(\frak D^*_w)&
=\left\langle\overline{I(W)} \ \Bigl | \ W:\text{words}, \ wt(W)=w \right\rangle_{\bold Q} \\
\intertext{by Proposition \ref{Explicit Formula}}
&=\left\langle\overline{Z(W)} \ \Bigl | \ W\in M :\text{words}
, \ wt(W)=w \right\rangle_{\bold Q} 
=NZ_w .
\end{align*}
This means that the linear map $\Psi_w$ is surjective.\par

{\it Preserving the depth filtration}:
We see easily that 
\[
\frak F^m\frak D^*_w=\left\langle W^*|_{\frak D_w}\Bigl |  \ W:\text{words}, \ wt(W)=w, \ dp(W)\leqslant m \right\rangle_{\bold Q}.
\]
Thus it also follows that 
\begin{align*}
\Psi_w(\frak F^m\frak D^*_w)&
=\left\langle\overline{I(W)} \ \Bigl |  \ W:\text{words}, \ wt(W)=w, \ dp(W)\leqslant m \right\rangle_{\bold Q} \\
\intertext{by Proposition \ref{Explicit Formula}}
&=\left\langle\overline{Z(W)} \ \Bigl | \ W\in M :\text{words}, \ wt(W)=w,  \ dp(W)\leqslant m \right\rangle_{\bold Q} \\
&=NZ_w^{\leqslant m}.
\end{align*}

This means that the linear map $\Psi_w$ strictly preserves 
the depth filtration.
\qed
\end{pf}
$\varPsi_{DR}$ stands for `de Rham' and `Drinfel'd'.

\subsection{Several corollaries and conjectures}\label{cor&conj}
The surjectivity of $\varPsi_{DR}$ implies 
\begin{cor}\label{Etingof}
\begin{enumerate}
\renewcommand{\labelenumi}{(\roman{enumi})}
\item   $dim_{\bold Q}NZ_w \leqslant dim_{\bold Q} \ \frak D_w $.
\item   More precisely, $dim_{\bold Q}NZ_w^{\leqslant m} 
\leqslant  dim_{\bold Q} \ \frak F^m\frak D^*_w 
=dim_{\bold Q} \ (\frak D_w /\frak F^{m+1}\frak D_w)  $.
\end{enumerate}
\end{cor}
M. Kaneko informed me that
Corollary \ref{Etingof}(i)
is also appeared in \cite{ES},
which was deduced from the action of the stable derivation algebra
on a certain torsor.
Proposition \ref{Takao} follows

\begin{cor}\label{TakaoZeta}
\begin{itemize}
\renewcommand{\labelenumi}{(\roman{enumi})}
\item $NZ_w^{\leqslant m}=NZ^{\leqslant m+1}_w$  if  $w\equiv m (mod 2)$.
\item  $NZ_w^{\leqslant m}=NZ_w$   if \qquad   $m>\frac{w}{2}-1$.
\item  $dim_{\bold Q} NZ_w^{\leqslant 1}\leqslant
\begin{cases}
1 & w=3,5,7,9,\dots    \\
0 & w:\text{otherwise}.
\end{cases} 
$
\item $dim_{\bold Q} NZ_w^{\leqslant 2}\leqslant
\begin{cases}
\  \ \ 0                & w:\text{odd}  \\  
 \ [\frac{w-2}{6}]  & w:\text{even}.
\end{cases} $
\end{itemize}
\end{cor}
Corollary \ref{TakaoZeta} (or at least a part of it)
seems to have been already found by Don Zagier (cf. \cite{Gon98}).
\begin{eg}\label{dim-fil}
From the table in Example \ref{Matsumoto}, we get

\begin{tabular}[t]{|c||c|c|c|c|c|c|c|c|c|c|c|c|}\hline
$w$ & 1&2&3&4&5&6&7&8&9&10&11&12\\ \hline
$dim_{\bold Q} \ NZ_w^{\leqslant 1}$ & 0 & 0 & 1 & 0 & ${\leqslant 1}$ & 0 & ${\leqslant 1}$ & $ 0 $ & ${\leqslant 1}$ & 0 & ${\leqslant 1}$  & 0\\ \hline
$dim_{\bold Q} \ NZ_w^{\leqslant 2}$ & 0 & 0 & 1 & 0 & ${\leqslant 1}$ & 0 & ${\leqslant 1}$ & ${\leqslant 1}$ & ${\leqslant 1}$ &  ${\leqslant 1}$ & ${\leqslant 1}$  &  ${\leqslant 1}$\\ \hline
$dim_{\bold Q} \ NZ_w^{\leqslant 3}$ & 0 & 0 & 1 & 0 & ${\leqslant 1}$ & 0 & ${\leqslant 1}$ & ${\leqslant 1}$ & ${\leqslant 1}$ &  ${\leqslant 1}$ & ${\leqslant 2}$  &  ${\leqslant 1}$\\ \hline
$dim_{\bold Q} \ NZ_w^{\leqslant 4}$ & 0 & 0 & 1 & 0 & ${\leqslant 1}$ & 0 & ${\leqslant 1}$ & ${\leqslant 1}$ & ${\leqslant 1}$ &  ${\leqslant 1}$ & ${\leqslant 2}$  &  ${\leqslant 2}$\\ \hline\hline
$dim_{\bold Q} \ NZ_w$& 0 & 0& 1 & 0 & ${\leqslant 1}$ & 0 & ${\leqslant 1}$ & ${\leqslant 1}$ & ${\leqslant 1}$ & ${\leqslant 1}$ & ${\leqslant 2}$ & ${\leqslant 2}$ \\ \hline
\end{tabular}
\end{eg}

\begin{eg}\label{11&12}
Using the computation table of the basis of $\frak D_w$ for $w=11$ and $12$ 
by H. Tsunogai,
one can extend the list appearing in \cite{GonECM}
to higher weights $11$ and $12$ as follows.
\begin{align*}
Z_{11}=&\langle
\pi^8\zeta(3),\pi^6\zeta(5),\pi^4\zeta(7),\pi^2\zeta(3)^3,\pi^2\zeta(9),
\zeta(3)^2\zeta(5),\zeta(3)\zeta(3,5), \\
&\qquad \qquad\qquad\qquad\qquad\qquad\qquad\qquad\qquad\qquad
\underset{dp=1}{\zeta(11)},
\underset{dp=3}{\zeta(2,1,8)}\rangle_{\bold Q}, \\
Z_{12}=&\langle
\pi^{12},\pi^{6}\zeta(3)^2,\pi^4\zeta(3)\zeta(5),\pi^4\zeta(3,5),
\pi^2\zeta(3)\zeta(7), \pi^2\zeta(5)^2,\pi^2\zeta(3,7),\\
&\qquad \qquad\qquad\qquad\zeta(3)^4,\zeta(3)\zeta(9),\zeta(5)\zeta(7),
\underset{dp=2}{\zeta(3,9)},
\underset{dp=4}{\zeta(2,1,1,8)}\rangle_{\bold Q}. \\
\end{align*}
\end{eg}

Let $\{ d'_k\}^\infty_{k=0}$ be the sequence determined by the following
series.
\[
\sum_{k=0}^\infty d'_k t^k:=\frac{1}{1-t^2}\prod_{w=1}^{\infty}
\frac{1}{(1-t^w)^{dim_{\bold Q}\frak D_w}} \ \ \ .
\]
Then $d'_k$ gives the dimension of the degree $k$-part of the graded polynomial algebra infinitely generated by $z_{i,j}$ 
$(1\leqslant i, \ 1\leqslant j \leqslant dim_{\bold Q}\frak D_i \ )$ with deg $z_{i,j}=i$ and 
$\pi^2$ with deg $\pi^2=2$ .
Thus we get the following dimension bounding of the MZV algebra.

\begin{cor}\label{dim-bd}
\[
dim _\bold Q Z_w \leqslant d'_w  \ \ \ \ \text{    for all } w.
\]
Here, the equality holds 
if and only if $\varPsi_{DR}:\frak D^*\centerdot\twoheadrightarrow NZ\centerdot$ is an isomorphism
and $Z \centerdot$ is polynomial algebra.
\end{cor}

\begin{pf}
It follows immediately from
\begin{align*}
\sum_{k=0}^\infty dim_{\bold Q}Z_k \cdot t^k & \leqslant
\frac{1}{1-t^2}\prod_{w=1}^{\infty}\frac{1}{(1-t^w)^{dim_{\bold Q}NZ_w}} \\
& \leqslant 
\frac{1}{1-t^2}\prod_{w=1}^{\infty}\frac{1}{(1-t^w)^{dim_{\bold Q}\frak D_w}}
=\sum_{k=0}^\infty d'_k t^k
\end{align*}
by Corollary \ref{Etingof}(i).
Here, for two formal power series $P(t)$, $Q(t)$ in $\bold Q[[t]]$,
$P(t)\geqslant Q(t)$ means that the formal power series $P(t)-Q(t)$ has all of 
its coefficients non-negative.
\qed
\end{pf}

From the table in Example \ref{Matsumoto}, we see that
$d'_w=d_w$ for $w\leqslant 12$.
Here, $d_w$ is the conjectured dimension of $Z_w$ in \S\S \ref{MZValg}.
Moreover, we find that

\begin{lem}
Suppose that Conjecture \ref{Deligne-Ihara}(3) for $\frak D\centerdot$ 
is true. Then 
\[
d'_w \leqslant d_w  \ \ \ \ \ \ \text{     for all } w.
\]
The equality holds if and only if Conjecture \ref{Deligne-Ihara}(2) for $\frak D\centerdot$ holds.
\end{lem}

\begin{pf}
Let $\Bbb F\centerdot=\underset{m\geqslant 1}{\oplus}\Bbb F_m $
be the free graded Lie algebra over $\bold Q$ generated by $e_m$ with 
deg $e_m=m$ $(m=3,5,7,\cdots)$.
From the assumption of this lemma, 
$dim_{\bold Q}\frak D_w \leqslant dim_{\bold Q}\Bbb F_w$
for all $w$.
So we get
\begin{align*}
\sum_{k=0}^\infty d'_k\cdot t^k 
& = \frac{1}{1-t^2}\prod_{w=1}^{\infty}
\frac{1}{(1-t^w)^{dim_{\bold Q}\frak D_w}} \\
& \leqslant \frac{1}{1-t^2}\prod_{w=1}^{\infty}
\frac{1}{(1-t^w)^{dim_{\bold Q}\Bbb F_w}} 
=\sum_{k=0}^\infty d_k\cdot t^k  \ \  .\\
\end{align*}
A short proof of the last equality was given in \cite{GonECM}.
The equality $d'_k=d_k$ for all $k$
holds if and only if
$dim_{\bold Q}\frak D_w=dim_{\bold Q}\Bbb F_w$ for all $w$,
which is equivalent to saying that
$\frak D\centerdot$ is the free graded Lie algebra generated by one element
in each degree $m (=3,5,7,\cdots)$.
\qed
\end{pf}

We can deduce the following proposition from Corollary \ref{dim-bd}.
\begin{prop}
Suppose that Conjecture \ref{Deligne-Ihara}(3) for $\frak D\centerdot$ is true, 
i.e. $\frak D\centerdot$ is generated by one element in each degree 
$m (=3,5,7,9,\cdots)$. Then
\[
dim_{\bold Q}Z_w \leqslant d_w  \ \ \ \ \ \ \text{       for all  } w.
\]
\end{prop}

Namely, the conjecture on the structure of the stable derivation algebra
which arose from the study of the Galois representation on 
$ \pi \sb 1(\bold P ^1\sb{\overline{\bold Q}}-\left\{0,1,\infty\right\})$
implies the upper-bounding part of the 
Dimension Conjecture of multiple zeta values (\S\S \ref{MZValg}) !\par
From Corollary \ref{dim-bd}, the validity of the Dimension Conjecture and the Conjecture \ref{Deligne-Ihara}(2),(3) for $\frak D\centerdot$ would imply that
$Z\centerdot$ might be a polynomial algebra (\cite{GonECM}) and
the following

\begin{conj}\label{deRham}
The surjective linear map
\[
\varPsi_{DR}:\Bigl(\frak D^*_\centerdot,\{\frak F^m\frak D^*_\centerdot\}_{m\in\bold N}\Bigr)
\twoheadrightarrow 
\Bigl(NZ\centerdot,\{NZ\centerdot^{\leqslant m}\}_{m\in\bold N}\Bigr)
\]
is an isomorphism as filtered graded $\bold Q$-vector space.
\end{conj}

But it may be hard to prove the injectivity,
since it implies deep results of transcendental number theory.

\begin{rem}
Particularly Conjecture \ref{deRham} implies that the dual vector space
$NZ^*\centerdot$ might form a sub-Lie algebra of $\frak D\centerdot$
by the Lie bracket $\delta$ of $\frak D\centerdot$.
This indicates that, for example, the map
$NZ^*_3\otimes NZ^*_5\to NZ^*_8$ induced from the Lie bracket
$\delta:\frak D_3\otimes \frak D_5\to \frak D_8$ might be an injection,
which means especially
$(0\leqslant)dim_{\bold Q}NZ_5\leqslant dim_{\bold Q}NZ_8(\leqslant 1)$ 
because $dim_{\bold Q}NZ_3=1$.
However showing this inequality looks a difficult problem 
(at least for the author) related to transcendental number theory.
\end{rem}

On the algebraic relations among MZV's, the author raises the following 
conjecture.
\begin{conj}\label{algebraic-relations}
All of the algebraic relations among the MZV's can be deduced from 
the relations of the Drinfel'd associator 
$(0)\sim(III)$ in Property II (\S\S \ref{Property II}).
\end{conj}
Suppose that Conjecture \ref{algebraic-relations} holds. 
Then Conjecture \ref{deRham} and Direct Sum Conjecture (\S\S \ref{MZValg})
also hold.
By the way, is it possible to deduce
all of the algebraic relations among the MZV's which were found 
till now from $(0)\sim(III)$?
This question does not seem so trivial at all.

\begin{rem}[(\cite{De}\S\S 18.13-18.17)]
P. Deligne deduced the Euler's formula
$\zeta(2n)=\frac{-(2\pi i)^{2n}}{2(2n)!}B_n$ 
($n\in \bold N$, $B_n$: the Bernoulli number)
from  $(0)\sim(II)$.
Since his proof is interesting, we give its brief sketch.\par
By $(I)$ and $(II)$, we get
\[
e^{\pi iC}\Phi_{KZ}(C,B)^{-1}e^{\pi iB}\Phi_{KZ}(A,B)e^{\pi iA}=
e^{-\pi iC}\Phi_{KZ}(C,B)^{-1}e^{-\pi iB}\Phi_{KZ}(A,B)e^{-\pi iA}
\]
He calculated the image of both hand sides of the above equality
modulo $exp \ \frak F^2\Bbb L^{\land}_{\Bbb C}$ in \cite{De} \S 18.16
as follows:
\begin{align*}
\text{(LHS)}\equiv & exp\Bigl[
\frac{e^{-\pi i(adA)}-1+e^{-\pi i(adA)}\{\pi i(adA)+2\sum_{n=1}^{\infty}\zeta(2n)(adA)^{2n}\}}
{ad A}(B)\Bigr] \\
&\mod{exp \ \frak F^2\Bbb L^{\land}_{\Bbb C}},\\
\text{(RHS)}\equiv & exp\Bigl[
\frac{e^{\pi i(adA)}-1+e^{\pi i(adA)}\{-\pi i(adA)+2\sum_{n=1}^{\infty}\zeta(2n)(adA)^{2n}\}}
{ad A}(B)\Bigr] \\
&\mod{exp \ \frak F^2\Bbb L^{\land}_{\Bbb C}},\\
\end{align*}
from which he obtained
\[
\sum\limits_{n=1}^\infty\zeta(2n)(adA)^{2n}(B)=
\sum\limits_{n=1}^\infty\frac{-(2\pi i)^{2n}}{2(2n)!}B_n(adA)^{2n}(B)\qquad .
\]
\end{rem}

\section{Some comparisons between `Galois Side' and `Hodge Side'}
For each $m\in\bold N$, the element $D=D_f\in \frak D_m$ determines
the unique rational number $c_m(D)$ by the congruence  
$f\equiv c_m(D)\cdot (ad x)^{m-1}(y) \quad mod  \ \frak F^2\Bbb L_m$.
Y. Ihara constructed 
a canonical $\bold Q$-linear map 
$c_m:\frak D_m\to \bold Q$
in \cite{Ih99}
which is defined by $D\mapsto c_m(D)$.
It is shown in \cite{Ih99} that $c_m\in\frak D^*_m$ 
is non-vanishing  
if and only if $m$ is odd and $m\geqslant 3$.\par

For $\sigma\in Gal({\overline{\bold Q}}/\bold Q)$, let 
$\tilde{\kappa}_m^{(l)}(\sigma)$ be the unique $l$-adic integer satisfying
\begin{equation*}
\prod_{a\in (\bold Z/l^n\bold Z)^\times}
\left(
\left[
\left\{ (\zeta_{l^n}^a-1)^{\langle a^{m-1}\rangle}\right\}^{\frac{1}{l^n}}
\right]^\sigma
\left/
\left\{(\zeta_{l^n}^{\chi(\sigma) a}-1)^{\langle a^{m-1}\rangle}\right\}
^{\frac{1}{l^n}}
\right.
\right)
=\zeta_{l^n}^{(l^{m-1}-1)\cdot\tilde{\kappa}_m^{(l)}(\sigma)}
\end{equation*}
for all $n\in\bold N$.
Here $\zeta_{l^n}=exp(\frac{2\pi i}{l^n})$ , $\langle a^{m-1}\rangle$ 
is the representative of $a^{m-1}\  mod \ l^n$ with $0<\langle a^{m-1}\rangle<l^n$ 
and $\chi$ is the $l$-adic cyclotomic character.
It is non-vanishing if and only if $m$ is odd and $m\geqslant 3$ (\cite{Sou}).
The map $\tilde{\kappa}_m^{(l)} \ (m\geqslant 3 , \text{ odd})$ is called the
{\it $m$-th Soul\'{e} element}.
It represents a non-trivial generator of 
$H^1_{\text{\'{e}t}}(Spec \ \bold Q,\bold Q_l(m))$,
which is of rank $1$ for odd $m\geqslant 1$ and rank $0$ for other $m$
$(\geqslant 1)$.\par

Taking the dual of the embedding
$\varPsi_l:\frak g^l\centerdot\hookrightarrow
\frak D\centerdot\otimes_\bold Q \bold Q_l$ 
which is associated with the Galois representation on 
the pro-$l$ fundamental group
$\pi_1^{l}(\bold P^1_{\overline{\bold Q}}-\left\{0,1,\infty\right\})$
(see \S\S \ref{connection}),
we get a surjective linear map
$\varPsi^*_l:\frak D^*_\centerdot\otimes_\bold Q \bold Q_l
\twoheadrightarrow\frak g^{l *}_\centerdot$ .
It is shown in \cite{Ih99} that
the image of $c_m\otimes_\bold Q  1$ by the map
$\varPsi^*_l:\frak D^*_\centerdot\otimes_\bold Q \bold Q_l
\twoheadrightarrow\frak g^{l *}_\centerdot$ \ 
is $\frac{1}{(m-1)!}\kappa_m^{(l)}$,
where $\kappa_m^{(l)}:\frak g^l_m\to\bold Q_l$ is
the induced $\bold Q_l$-linear map by $\tilde{\kappa}_m^{(l)}$.
On the other hand, we find that
the image of $c_m$ by
$\varPsi_{DR}:\frak D^*_\centerdot\twoheadrightarrow NZ\centerdot$ \ 
is $\overline{-\zeta(m)}$ by Proposition \ref{Explicit Formula}(a).
We can make a comparison between `the de Rham world' and
`the $l$-adic world' as follows.

\begin{align*}
&\boxed{\text{de Rham}} \qquad \qquad \qquad\qquad \qquad\quad
\boxed{\text{$l$-adic}} \\
&\qquad\quad  
NZ\centerdot \
\overset{\varPsi_{DR}}{\twoheadleftarrow}
\quad
\frak D^*_\centerdot  \\
&\qquad\qquad\qquad \qquad
\frak D^*_\centerdot\otimes_\bold Q \bold Q_l \ \ 
\overset{\varPsi^*_l}{\twoheadrightarrow} \qquad
\frak g^{l*}_\centerdot \\
&\qquad \overline{-\zeta(m)}\ 
\longleftarrow\!\!\mbox{\rule[0.9pt]{0.4pt}{4.1pt}} 
\qquad
c_m \qquad \
\longmapsto 
\frac{1}{(m-1)!}\kappa_m^{(l)}
\end{align*}

Conjecture \ref{deRham} and Conjecture \ref{Deligne-Ihara} expect that both
$\varPsi_{DR}$ and $\varPsi^*_l$ are isomorphisms.


\end{document}